\documentclass{jktrX}

\def\draftnote{ }
\def\today{}
\arraycolsep=\labelsep

\usepackage{graphicx}

\newcommand{\bysame}{\mbox{\rule{3em}{.4pt}}\,}

\usepackage{amssymb}
\usepackage{amsmath}
\usepackage{latexsym}

\newcommand\ZZ{{\Bbb Z}}

\newcommand\QQ{{\Bbb Q}}
\newcommand\qed{\hfill $\Box$ \hfill \\} 
\renewcommand\qed{\hfill $\Box$ \hfill \\}

\def\marusen{\unitlength.1em
  \begin{minipage}{10\unitlength}
    \begin{picture}(10,10)
      \put(5,6){\circle{6}} %\put(5,6){\circle{3}}
      \qbezier(3,8)(7,4)(7,4)
    \end{picture}
  \end{minipage}
}
\def\marubackslash{{\small \marusen}}

\newtheorem{thm}{Theorem}[section]

\newtheorem{prop}[thm]{Proposition}

\newtheorem{yosou}[thm]{Conjecture}
 \newtheorem{dfn}[thm]{Definition}

 \newtheorem{ex}[thm]{Example}

\title{Twisted Alexander polynomials of 
2-bridge knots associated to
metabelian representations}

\author{Mikami Hirasawa}
\address{Department of Mathematics,
Nagoya Institute of Technology,\\
Nagoya Aichi 466-8555 Japan\\
{\it E-mail: hirasawa.mikami@nitech.ac.jp}
}

\author{Kunio Murasugi}
\address{Department of Mathematics,
University of Toronto,\\
Toronto, ON M5S2E4 Canada\\
{\it E-mail: murasugi@math.toronto.edu}
}

\begin{document}

\maketitle
 \pagestyle{myheadings}
 \markboth{
Twisted Alexander polynomials 
associated to  metabelian representations
}
 {M. Hirasawa \& K. Murasugi }

\begin{abstract}
Suppose the knot group $G(K)$ of a knot $K$ has a non-abelian
representation $\rho$ on $A_4 \subset GL(4,\ZZ)$. 
We conjecture that the twisted
Alexander polynomial of $K$ associated to $\rho$ is of the form:
$\displaystyle{\left[\dfrac{\Delta_K(t)}{1-t}\right]\varphi(t^3)}$,
where $\varphi(t^3)$ is an integer polynomial in $t^3$. 
We prove the conjecture for $2$-bridge knots $K$ 
whose group $G(K)$ can be mapped onto 
a free product $\ZZ/2*\ZZ/3$.
Later, we discuss more general metabelian representations 
of the knot groups and propose a similar conjecture on the 
form of the twisted Alexander polynomials.
\end{abstract}


\section{Introduction}
In this paper, which is a sequel of \cite{HM1}, \cite{HM2},
we consider non-abelian representations of the knot group
$G(K)$ of a knot $K$ on metabelian groups.
First we study a representation of $2$-bridge knot groups
on the alternating group $A_4$ of order $12$, 
the simplest metabelian group, which we call
an $A_4$-representation.
It is shown in \cite{R} and \cite{Ha} %[R] and [2] 
that $G(K)$ has a $A_4$-representation if and
only if
$\Delta_K(\omega) \Delta_K (\omega^2) \equiv 0$ ($mod$ 2), 
where $\Delta_K(t)$
is the Alexander polynomial of $K$ and 
$\omega$ is a primitive cubic root~of~$1$.
For a $2$-bridge knot $K(r)$, Heusner gives a nice criterion for $G(K(r))$
to have
an $A_4$-representation in terms of the degree of $\Delta_K(t)$ \cite{He}.
%[H]. In
In the later section,
we discuss briefly the general metabelian representation of the knot groups.

Now, let $\rho: G(K) \rightarrow A_4$ be an $A_4$-representation of
$G(K)$.
We may assume that one meridian generator maps to $(123)$. 
Let $K(r)$, $r \in \QQ$, $0 < r < 1$,
be a $2$-bridge knot with a Wirtinger presentation
%(1.1) $G(K(r))$ = $<x, y: R>$, $R = W x W^{-1} y^{-1}$.
\begin{equation} %1.1
G(K(r))=\langle x, y|R\rangle, R=WxW^{-1}y^{-1}.
\end{equation}
Suppose that $G(K(r))$ has an $A_4$-representation $\rho$ such that
%(1.2) $\rho (x)$ = (123) and $\rho (y)$ = (142).
\begin{equation}%1.2
\rho(x)=(123)\ {\rm and}\ \rho(y)=(142).
\end{equation}
Let $\xi:A_4 \rightarrow GL(4,\ZZ)$ be the permutation
representation of $A_4$. 
Then $\xi$ is equivalent to $\widetilde{\xi}$:
\begin{equation}%1.3
\widetilde{\xi}(123)=
\left(\begin{array}{rrrr}
1&-1&0&0\\
0&-1&1&0\\
0&-1&0&0\\
0&-1&0&1
\end{array}
\right)
\ {\rm and}\ 
\widetilde{\xi}(142)=
\left(
\begin{array}{rrrr}
1&0&0&-1\\
0&0&0&-1\\
0&0&1&-1\\
0&1&0&-1
\end{array}
\right)
\end{equation}
and hence, 
$\xi_0:
A_4 \longrightarrow  GL(3,\ZZ)$ 
given by
\begin{equation}%1.4
\xi_0(123)=\left(
\begin{array}{rrr}
-1&1&0\\
-1&0&0\\
-1&0&1
\end{array}\right)
\ {\rm and}\
\xi_0(142)=
\left(
\begin{array}{rrr}
0&0&-1\\
0&1&-1\\
1&0&-1
\end{array}\right)
\end{equation}
defines an irreducible representation of $A_4$. 
Combining them, we have a
$A_4$-representation of $G(K(r))$, 
$\rho_0=\rho \circ \xi_0 : 
G(K(r)) \longrightarrow GL(3, \ZZ)$, and the twisted Alexander polynomial
$\widetilde{\Delta}_{\rho_0, K(r)} (t)$ of $K(r)$ associated to $\rho_0$
is defined. 

From the forms of $\widetilde{\xi}$ and $\xi$,
we see immediately that the twisted Alexander polynomial
$\widetilde{\Delta}_{\rho,K}(t)$ is the product of
$\widetilde{\Delta}_{\iota,K}(t)$ and $\widetilde{\Delta}_{\rho_0, K}(t)$, where $\iota$
is a trivial representation. Since $\widetilde{\Delta}_{\iota,K}(t)=\Delta_{K}(t)/(1-t)$, the conjecture stated in the abstract is rephrased as follows:

%
%We propose the following conjecture on the form of these
%twisted Alexander polynomials. 
%(For a more general conjecture, see section 6.)

%%Conjecture A. 
\medskip
\noindent
{\bf Conjecture A.}
{\it
For a 2-bridge knot $K(r)$ with an $A_4$-representation,
$\widetilde{\Delta}_{\rho_0, K(r)} (t)$ 
is an integer polynomial in $t^3$
(up to $\pm t^k$ ), namely,
\begin{equation}%1.5
\widetilde{\Delta}_{\rho_0, K(r)} (t) =
\pm t^k  \varphi(t^3),\ 
{\it for\ some\ integer\ polynomial}\ \varphi(t).
\end{equation}
}%end \it

We prove this conjecture for $K(r)$ in $H(3)$, where $H(3)$ is the set of
$K(r)$ such that $G(K(r))$ maps onto a non-trivial free
product $\ZZ/2 * \ZZ/3$. 
The proof will be given in Section 2 through Section 5.
Since our proof is similar to those given in the previous two papers
\cite{HM1} and
\cite{HM2}, 
we will skip some details in our argument.

In the last Section 6, we discuss general metabelian representations and
state a similar conjecture on the form of the twisted Alexander polynomials.
We give several examples that justify our conjecture.

For convenience, we draw a diagram below consisting of various groups and
homomorphisms.

\begin{center}
$
\begin{array}{crccc}
& & & & GL(4,\ZZ)\\
& & &\nearrow \mbox{\large $\xi$} &\\
G(K)& \underset{\mbox{\large $\rho$}}{\longrightarrow}& A_4&\underset{\mbox{\large $\xi_0$}}{\longrightarrow}&GL(3,\ZZ)\\
& &\downarrow& &\\
& & \ZZ A_4&\underset{\mbox{\large $\widetilde{\xi}_0$}}{\longrightarrow}&M_{3,3}(\ZZ)\\
 & & & & \\
\ZZ F(x,y)&\underset{\mbox{\large $\widetilde{\rho}$}}{\longrightarrow}&
\widetilde{A}(x,y)& &\\
\downarrow\nu^*& & & &\\
\mbox{$[\ZZ F(x,y)][t^{\pm 1}]$}&
\underset{\mbox{\large $\rho^*$}}{\longrightarrow}&
\widetilde{A}(x,y)[t^{\pm 1}]&
\underset{\mbox{\large $\eta$}}{\longrightarrow}&
GL(3,\ZZ [t^{\pm 1}])\\
\end{array}
$
\end{center}
Here, 
(1) $\rho_0=\rho \circ \xi_0$,
(2) $\nu^*(\prod_{i=1}^m x^{k_i}
y^{\ell_i})=
(\prod_{i=1}^m x^{k_i} y^{\ell_i}) t^q$, 
$q=\sum_{i=1}^m
k_i + \sum_{i=1}^m \ell_i$.

%Sec2.
\section{$\ZZ$-algebra $\widetilde{A}(x,y)$}

Denote $x=(123)$ and $y=(142)$. 
Let $X= \xi_0 (x)$ and $Y= \xi_0 (y)$.
We define an algebra 
$\widetilde{A}(x,y)$ using the group algebra 
$\ZZ A_4$ as follows.

Let  $\widetilde{\xi}_0: 
\ZZ A_4 \longrightarrow M_{3,3}(\ZZ)$ 
be a linear extension of $\xi_0$ and
$S =(\widetilde{\xi}_0)^{-1}(0)$ be the kernel of 
$\widetilde{\xi}_0$.
Then $\ZZ A_4/S$ is a non-commutative
$\ZZ$-algebra, denoted by $\widetilde{A}(x,y)$.

%Proposition 2.1. 
\begin{prop}\label{prop:2.1}
In $\widetilde{A}(x,y)$, the following formulas hold.
\begin{align}%2.2
&(1)\ x^3 = y^3 =(xy)^3 = 1,\nonumber\\
&(2)\ xyx = yxy,\nonumber\\
&(3)\ (xy^{-1})^2 = 1,\nonumber\\
&(4)\ xyx = x^{-1} y^{-1} x^{-1},\nonumber\\
&(5)\ (x+y)^2 = (x^{-1} + y^{-1})^2 = 0,\nonumber\\
&(6)\ xyx(x+y) = (x+y)xyx = -(x+y),\nonumber\\
&(7)\ xyx(x^{-1} + y^{-1}) = 
(x^{-1} + y^{-1})xyx = - (x^{-1} + y^{-1}),\nonumber\\
&(8)\ (x +y)(x^{-1} + y^{-1}) + (x^{-1} + y^{-1})(x+y) 
= 2(1-xyx),\nonumber\\
&(9)\ xy + yx = -(x^{-1} + y^{-1})\  {\rm and}\ 
x^{-1} y^{-1} + y^{-1} x^{-1} =- (x+y).
\end{align}
\end{prop}

Only a few formulas below are needed to prove Proposition \ref{prop:2.1}, and details are omitted. \qed
%We need only few formulas below to prove Proposition \ref{prop:2.1}
%Proposition 2.1.

\arraycolsep=0.14em

\begin{equation}%2.2
X + Y = \left[\begin{array}{rrr}
-1&1&-1\\
-1&1&-1\\
0&0&0
\end{array}
\right],
X^{-1} + Y^{-1}= \left[\begin{array}{rrr}
-1&-1&1\\
0&0&0\\
-1&-1&1
\end{array}
\right]\ {\rm and}\ 
XYX=\left[
\begin{array}{rrr}
-1&0&0\\
-1&0&1\\
-1&1&0
\end{array}
\right]
\end{equation}

\arraycolsep=0.5em

%
%
%%%\begin{align}%2.2
%%%X + Y  &= \left(\begin{array}{rrr}
%%%-1&1&-1\\
%%%-1&1&-1\\
%%%0&0&0
%%%\end{array}
%%%\right),\  \nonumber\\
%%%X^{-1} + Y^{-1} &= \left(\begin{array}{rrr}
%%%-1&-1&1\\
%%%0&0&0\\
%%%-1&-1&1
%%%\end{array}
%%%\right)\ {\rm and} \nonumber\\
%%%XYX&=\left(
%%%\begin{array}{rrr}
%%%-1&0&0\\
%%%-1&0&1\\
%%%-1&1&0
%%%\end{array}
%%%\right)
%%%\end{align}
%%%%(2.2)    $X + Y$  = $[-1,1,-1 \mid -1,1,-1 \mid 0,0,0]$,
%%%%
%%%%         $X^{-1} + Y^{-1}$ =  $[-1,-1,1 \mid 0,0,0 \mid -1,-1,1]$ and
%%%%
%%%%         $XYX$ = $[-1,0,0 \mid -1,0,1 \mid -1,1,0]$.
%%%Details are omitted.\qed

Now, let $L(t)$ be the set of Laurent polynomials in 
$t$ with coefficients in $\widetilde{A}(x,y)$.
Write  $f(t)= {\displaystyle \sum_{- \infty < j < \infty}}
 d_j
t^j,\ d_j \in \widetilde{A} (x,y)$.

%Definition 2.2.  
\begin{dfn}
A polynomial $f(t)$ in $L(t)$ is called {\it twin} if $f(t)$
satisfies the following conditions:
%(2.3) 
\begin{align}
&(1)\  {\rm If}\ 
j \equiv  0\ ({\rm mod}\  3),\ 
{\rm then}\ d_j = c_j + c_j^{\prime}
xyx,\  {\rm where}\ c_j , c_j^{\prime} \in \ZZ.\nonumber\\
&(2)\  {\rm If}\ 
j \equiv 1\ ({\rm mod}\  3),\
{\rm then}\ 
d_j = a_j (x+y),\ a_j  \in \ZZ.\nonumber\\
&(3)\   {\rm If}\ 
j \equiv 2\ ({\rm mod}\ 3),\ {\rm then}\ 
d_j = b_j (x^{-1} + y^{-1}),\ b_j\in \ZZ.
\ {\rm Furthermore},\nonumber\\
&(4)\   {\rm for\ any}\ j,  
a_{3j+1} = b_{3j+2}.
\end{align}
The set of twin polynomials is denoted by $T(t)$.
\end{dfn}

%Proposition 2.3.   
\begin{prop}\label{prop:2.3}
$T(t)$ is a non-commutative subring of $L(t)$.
\end{prop}

\noindent
{\it Proof.}  Let $f(t)$ and $g(t)$ be twin polynomials. 
Obviously, $f(t) \pm g(t)$ is twin.
To show $f(t)g(t)$ is twin, 
it is enough to show that for any $m \in \ZZ$,
%(2.4) 
\begin{align}
&(1)\ f(t) t^{3m} \in  T(t)\ {\rm and}\ 
f(t)xyx t^{3m}  \in T(t).\nonumber\\
&(2)\  \left\{(x+y) t^{3m+1} + (x^{-1} + y^{-1}) t^{3m+2}
\right\}\left\{(x+y)
t^{3\ell+1} + (x^{-1} + y^{-1}) t^{3\ell+2}\right\} \in T(t).
\end{align}
These formulas follow from (2.1) (5)-(8).
\qed

%3. Main Theorem
\section{Main Theorem.}
Let $K(r)$ be a 2-bridge knot in $H(3)$.  
Then the  continued fraction of
$r =
\beta / \alpha,\ \alpha \equiv \beta \equiv 1$ 
(mod $2$), is of the following form. (See \cite{GR}, \cite{ORS}) %[G,ORS]:

\begin{equation*}
r =[3k_1,2m_1, 3k_2, 2m_2, \cdots, 2m_{q-1}, 3k_q].
\end{equation*}

Consider a Wirtinger presentation of $G(K(r))$:
\begin{equation*}
G(K(r)) = \langle x,y |  R \rangle,   
R = W x W^{-1} y^{-1}.
\end{equation*}
\noindent
Since $K(r) \in  H(3)$, $R$ is a product of conjugates of $R_0
=xyxy^{-1}x^{-1}y^{-1}$, namely

%(3.1)        
\begin{equation}
R = \prod_{j=1}^{m}u_j R_0^{\epsilon_j} u_j^{-1},
\end{equation}
where $u_j$ are words in the free group $F(x,y)$ freely generated by $x$
and $y$,  and $\epsilon_j = \pm 1$.
Therefore,  
\begin{equation*}
\frac{\partial R}{\partial x}
=\sum_{j=1}^{m} \epsilon_j u_j
\frac{\partial{R_0}}{\partial x},
\ {\rm where}\ 
\frac{\partial}{\partial x}\ 
{\rm denotes\ the\ free\ derivative}.
\end{equation*}
Let  $\widetilde{\rho}: 
\ZZ F(x,y) \longrightarrow \widetilde{A} (x,y)$ be
an algebra homomorphism defined by
\begin{equation*}
\widetilde{\rho} (x)=x\ {\rm and}\ 
\widetilde{\rho} (y) = y.
\end{equation*}
We write $\lambda (r) = 
\widetilde{\rho} (\sum_{j} \epsilon_j u_j) \in
\widetilde{A} (x,y)$ and $\lambda^* (r)
=(\nu^* \circ \rho^*)(\lambda
(r)) \in \widetilde{A} (x,y)[t^{\pm 1}]$.
Then the following is shown in \cite{HM1}: %[HM 1]

%(3.2) 
\begin{equation}
\widetilde{\Delta}_{\rho_0,K(r)} (t) 
=\left\{\det \widetilde{\xi}_0(\lambda^*(r))\right\}
\widetilde{\Delta}_{\rho_0,K(1/3)} (t).
\end{equation}
We note that 
$\widetilde{\Delta}_{\rho_0,K(1/3)} (t)$ = $1-t^3$.

Now our main theorem is the following:

%Theorem 3.1.   
\begin{thm}\label{thm:3.1}
If $K(r) \in H(3)$, then
$\widetilde{\Delta}_{\rho_0,K(r)} (t) 
= \varphi (t^{\pm 3})$,
where $\varphi (t)$ is an integer polynomial.
\end{thm}

Since 
$\widetilde{\Delta}_{\rho_0,K(1/3)} (t)=1-t^3$, 
the theorem is a consequence of 
Proposition \ref{prop:3.2} below:

%Proposition 3.2.  
\begin{prop}\label{prop:3.2}
Under the same conditions of Theorem \ref{thm:3.1}, 
we have; %3.1,
\begin{equation*}
\det[ \widetilde{\xi}_0 (\lambda^* (r))] = \varphi_0 (t^{\pm 3}),\
{\it for\ some\ integer\ polynomial}\ \varphi_0 (t).
\end{equation*}
\end{prop}

%4.
\section{Proof of Proposition 3.2. (I)}

In this section, we prove some basic formulas needed 
to prove the main theorem.  However,
since these formulas are mostly technical, 
we omit some details. (See \cite{HM2}.)%[HM2].)

Now, we denote;
%
%(4.1)   
\begin{align}
&(1)\ Q_0 (t) = 1,\nonumber\\
&(2)\ {\rm For}\  m \geq 1,\nonumber\\
&\ \ \ 
(a)\  Q_m (t) = 1 +(yx) t^2 + (yx)^2 t^4 
+ \cdots + (yx)^m t^{2m},\nonumber\\
&\ \ \ 
(b)\ Q_{-m} (t)
= (yx)^{-m} t^{-2m} Q_{m-1}(t) =
(x^{-1} y^{-1}) t^{-2} + (x^{-1} y^{-1})^2 t^{-4} 
+ \nonumber\\
&\ \ \ \ \ \ \ \ \cdots + (x^{-1}
y^{-1})^m t^{-2m}.
\end{align}

  We claim first the following:

%Proposition 4.1.  
\begin{prop}\label{prop:4.1}
Let $k \geq 0$, then we have;
%(4.2)  
\begin{align}
&(1)\  
y^{-1} t^{-1}[(1-yt) Q_{3k+1}(t) yt + (yx)^{3k+2}
t^{6k+4}](1-xt) \in T(t).\nonumber\\
&(2)\  
y^{-1}t^{-1}(1-yt) Q_{3k+2}(t) yt (1-xt) \in T(t).\nonumber\\
&(3)\  
y^{-1} t^{-1} [(1-yt) Q_{-(3k+1)}(t) yt - (x^{-1}
y^{-1})^{3k+1} t^{-(6k+2)}](1-xt)  \in T(t).
\nonumber\\
&(4)\  y^{-1} t^{-1} (1-yt) Q_{-(3k+3)}(t) yt (1-xt)  
\in T(t).
\end{align}
\end{prop}

\noindent
{\it Proof.}  For simplicity, we denote  
${\bf a}  = x+y$ and ${\bf b} = x^{-1} +y^{-1}$.
A proof will be done by induction on $k$.  First, straightforward
computations prove for the initial cases.
\begin{align*}
(1)& \
\ y^{-1} t^{-1}\{(1-yt) Q_1(t) yt + (yx)^2 t^4\}(1-xt)\\
&=1 -{\bf a}
t - {\bf b} t^2 - xyx t^3 \in T(t).\\
(2)&  \
\ y^{-1} t^{-1} (1-yt) Q_2 (t) yt (1-xt)\\
&= 1 -{\bf  a} t - {\bf b}
t^2 +2 xyx t^3 - {\bf a} t^4
- {\bf b} t^5 + t^6 \in T(t).\\
(3)&  \ 
\ y^{-1} t^{-1} \{(1-yt) Q_{-1} (t) yt 
- (x^{-1} y^{-1}) t^{-2}\}(1-xt)\\
&=- xyx t^{-3} - {\bf a} t^{-2}
- {\bf b} t^{-1} + 1 \in T(t).\\
(4)&\  
\  y^{-1} t^{-1} (1-yt) Q_{-3} (t) yt (1-xt)\\
&= t^{-6} - {\bf a}
t^{-5} - {\bf b} t^{-4} - 2xyx t^{-3}
- {\bf a} t^{-2} - {\bf b} t^{-1} +1 \in  T(t).
\end{align*}

Now suppose the formula hold for 
$k = k$ and we prove them for $k=k+1$.

{\it Proof of (1)}. 
Since $Q_{3k+4}(t) = Q_{3k+1} + (yx)^{3k+2} t^{6k+4}
Q_2(t)$, we see
\begin{align*}
&\ \ \ \ y^{-1} t^{-1} \{(1-yt) Q_{3k+4}(t) yt + 
(yx)^{3k+5} t^{6k+10}\}(1-xt)\\
&= 
y^{-1} t^{-1} \{(1-yt) Q_{3k+1} (t) yt + 
(yx)^{3k+2} t^{6k+4} + (1-yt)
(yx)^{3k+2} t^{6k+4}
Q_2 (t) yt     \\
&\ \ \ \ -(yx)^{3k+2} t^{6k+4} +(yx)^{3k+5} t^{6k+10}
\}(1-xt).
\end{align*}
Therefore, by induction hypothesis, it suffices to show that
%(4.3)    
\begin{align}
&y^{-1} t^{-1}\{(1-yt) (yx)^{3k+2} t^{6k+4} Q_2 (t) 
yt  \nonumber\\
&\ \ \ -(yx)^{3k+2} t^{6k+4} 
+(yx)^{3k+5} t^{6k+10}\}(1-xt) \in T(t),
\end{align}
or equivalently

%(4.4)
\begin{align}
y^{-1}\{(1-yt) (yx)^2 Q_2 (t) yt  
-  (yx)^2 + (yx)^2 t^6 \}(1-xt)
\in T(t).
\end{align}
However, it is straightforward to show that
\begin{align*}
&\{y^{-1} (1-yt) (yx)^2 Q_2 (t) yt - xyx  + 
xyx t^6\}(1-xt)\\
= &
-xyx - {\bf a} t - {\bf b}
 t^2 +2 t^3 - {\bf a} t^4 - {\bf b} t^5 -  xyx
t^6 \in T(t).
\end{align*}
Since other cases can be handled in the same way, details will be omitted.
\qed

Now, to prove Proposition 3.2, we need explicit recursion formulas for
$\lambda (r)$.
 To obtain such formulas, we write
 \begin{equation*}
r_q=[3k_1, 2m_1,3k_2, 2m_2, \dots, 2m_{q-1}, 3k_q].
\end{equation*}
Then  $\lambda (r_q)$ is exactly $w_{2q-1}$ in \cite{HM1}.
%[HM 1]. 
Using formula there, we can compute $\lambda (r_q)$ inductively.  
We only state formulas without proof.

%Proposition 4.2.   
\begin{prop}\label{prop:4.2}
Let
$r_q = [3k_1, 2m_1,3k_2, 2m_2, \dots, 2m_{q-1},3k_q]$.

Case 1.  $k_q > 0$.

(1) If $k_q = 2s$, $s \geq 1$, then

\begin{align*}
\lambda (r_q)=
&(1-y) Q_{3s-1} y 
  \{\sum_{j=1}^{q-1} m_j (x-1) y^{-1}
\lambda (r_j)\} \\
&+ (yx)^{3s} \lambda (r_{q-1}) -
\sum_{j=1}^s (yx)^{3s-3j+2}  
+ \sum_{j=1}^s (yx)^{3s-3j} y.
\end{align*}

(2) If $k_q = 2s -1$, $s \geq 1$, then

 \begin{align*}
 \lambda (r_q)=
 &\{(1-y) Q_{3s-2} y  + (yx)^{3s-1}\} 
 \sum_{j=1}^{q-1} m_j
(x-1) y^{-1} \lambda (r_j) \\
&-  (yx)^{3s-1} y^{-1} \lambda (r_{q-1})
+ \sum_{j=1}^s (yx)^{3s-3j} y  
-  \sum_{j=1}^{s-1} (yx)^{3s-3j-1}.
\end{align*}

Case 2.  $k_q < 0$.

(1)  If $k_q$ = $- 2s$, $s \geq 1$, then

\begin{align*}
\lambda (r_q)
 =&
  (y-1) Q_{-3s} y \sum_{j=1}^{q-1} m_j (x-1) y^{-1}
\lambda (r_j) + (x^{-1} y^{-1})^{3s} 
\lambda (r_{q-1})\\
& -\sum_{j=1}^s (x^{-1} 
y^{-1})^{3s-3j+2} x^{-1}  +  \sum_{j=1}^s
(x^{-1} y^{-1})^{3s-3j+1}.
\end{align*}

(2) If $k_q = -(2s +1)$, $s \geq 0$, then

 \begin{align*}
 \lambda(r_q)=&
 \{(y-1) Q_{-(3s+1)}y + (x^{-1} y^{-1})^{3s+1}\}
\sum_{j=1}^{q-1} m_j (x-1) y^{-1} \lambda (r_j)\\ 
&
-(x^{-1}y^{-1})^{3s+1} y^{-1} \lambda (r_{q-1}) 
+ \sum_{j=0}^s (x^{-1}
y^{-1})^{3s-3j+1} \\
&- \sum_{j=0}^s 
(x^{-1} y^{-1})^{3s-3j+2} x^{-1}.
\end{align*}

\end{prop}

Using Proposition \ref{prop:4.2}, %4.2, 
we can show the following key proposition.

%Proposition 4.3.  
\begin{prop}\label{prop:4.3}
For any $r_q$, $q \geq 1$,
$y^{-1} t^{-1}\lambda^*(r_q)$ is twin.
\end{prop}

\noindent
{\it Proof.} 
Since  other cases can be proven in the same way, 
we prove only one case: $k_q = 2s$, $s \geq 1$.
We use induction argument.

Let $q = 1$. Then $r=[6s]$, $s \geq 1$, and
\begin{equation*}
\lambda (r_1)=
- \sum_{j=1}^s (yx)^{3s-3j+2} + \sum_{j=1}^s
(yx)^{3s-3j} y.
\end{equation*}
Therefore
\begin{align*}
y^{-1}t^{-1}\lambda^*(r_1)&=
y^{-1} t^{-1}\{ - \sum_{j=1}^s (yx)^2
t^{6s-6j+4} + \sum_{j=1}^s y t^{6s-6j+1}\}\\
&=- \sum_{j=1}^s xyx t^{6s-6j+3} + \sum_{j=1}^s t^{6s-6j},
\end{align*}
that is obviously twin.\\
Now suppose  $y^{-1} t^{-1} \lambda^* (r_j)$ is twin for
$j=1,2,\dots,q-1$. Then
\begin{align*}
y^{-1} t^{-1} \lambda^* (r_q)
&=y^{-1} t^{-1} (1-yt) Q_{3s-1}(t) yt \left\{
\sum_{j=1}^{q-1} m_j (xt-1) y^{-1} t^{-1} \lambda^* (r_j)\right\}\\
&\ \ \ + y^{-1} t^{-1} (yx)^{3s} t^{6s} \lambda^*(r_{q-1}) - 
 y^{-1} t^{-1}
\sum_{j=1}^s (yx)^{3s-3j+2} t^{6s-6j+4}\\ 
&\ \ \ +y^{-1} t^{-1} \sum_{j=1}^s (yx)^{3s-3j} y t^{6s-6j+1}\\
&=\sum_{j=1}^{q-1} m_j \{ y^{-1} t^{-1} (1-yt) Q_{3s-1}(t) yt (xt-1)\}
y^{-1} t^{-1} \lambda^*(r_j)\\
&\ \ \ 
+ t^{6s} \{y^{-1} t^{-1} \lambda^*(r_{q-1})\} - 
\sum_{j=1}^s xyx
t^{6s-6j+3} + \sum_{j=1}^s t^{6s-6j}.
\end{align*}
Since both $y^{-1} t^{-1} \lambda^* (r_j)$ and $y^{-1} t^{-1} (1-yt)
Q_{3s-1} (t) yt (xt-1)$ are twin by
(4.2)(2), Proposition \ref{prop:4.3} %4.3 
follows for this case.      \qed

%5.
\section{Proof of Proposition 3.2.(II)}

To show $\det[\widetilde{\xi}_0(\lambda^*(r))]
= \varphi_0(t^{\pm 3})$, it suffice to
show $\det[\widetilde{\xi}_0(y^{-1}t^{-1} \lambda^*(r))]
=\varphi_1(t^{\pm 3})$.
Now, since $y^{-1}t^{-1} \lambda^*(r)$ is twin, we can write
\begin{align*}
y^{-1}t^{-1} \lambda^*(r)&=\ 
\sum_{-\infty < j <\infty}(c_j +
c_j^{\prime} xyx ) t^{3j} \\
&\ \ + \sum_{-\infty < j < \infty} a_j(x+y)
t^{3j+1}\\
&\ \ + 
\sum_{-\infty < j < \infty} b_j(x^{-1} + y^{-1}) t^{3j+2},
\end{align*}
where $a_j, b_j, c_j$ and $c_j^{\prime}$ are integers and $a_k$= $b_k$ for
all  $k$.\\
Denote $A=\sum_{j} a_j t^{3j+1}$, 
$C=\sum_{j}c_j t^{3j}$ and
$C^{\prime}=\sum_{j}c_j^{\prime} t^{3j}$.
Then, 
\begin{equation*}
y^{-1} t^{-1} \lambda^*(r)=
C +C^{\prime} xyx +A(x+y) + A(x^{-1}
+ y^{-1})t.
\end{equation*}

Using (2.2), we see that
\begin{align*}
D=&\det[ \widetilde{\xi}_0 (y^{-1} t^{-1} \lambda (r))]\\
=&\det\left(\begin{array}{ccc}
-A - At+C - C^{\prime}&\ A-At&\ -A+At \\
  -A - C^{\prime}&\ A+C&\ -A+C^{\prime} \\
 -At -C^{\prime}&\ -At+C^{\prime}&\ At+C
 \end{array}\right).
 \end{align*}
First add the second column to the third, 
then subtract the second row from
the third, and we have:
\begin{align*}
D&
=(C+C^{\prime}) \{(-A -At+C -C^{\prime})
(-A-At -C +C^{\prime}) -
(A-At)^2\}\\
&=(C +C^{\prime}) \{4A^2 t - (C-C^{\prime})^2\}.
\end{align*}
Since $A=t \sum_{j} a_j t^{3j}$, 
we see $A^2=t^2 \{\sum_{j} a_jt^{3j}\}^2$ and thus
\begin{equation*}
D=\left\{ \sum_{j} (c_j + c_j^{\prime}) t^{3j}\right\}
\left\{4t^3( \sum_{j} a_j t^{3j})^2 -  ( \sum_{j}(c_j
 - c_j^{\prime})t^{3j})^2\right\}.
 \end{equation*}
Therefore $D$ is a polynomial in $t^{\pm 3}$. \\
This proves Proposition \ref{prop:3.2}, %3.2
and the proof of Theorem \ref{thm:3.1} %3.1 
is now complete.

\begin{ex}
The following examples justify our main theorem.\\
(1) For  $r=1/9,\widetilde{\Delta}_{\rho_0,K(r)} (t)
=(1-t^3)(1-t^3 + t^6 )(1 + t^3 + t^6)^2$.\\
(2) For $r = 5/27, 
\widetilde{\Delta}_{\rho_0,K(r)} (t) 
=(1-t^3)(4+7t^3 +4t^6)$.\\
(3) For  $r = 7/39,
\widetilde{\Delta}_{ \rho_0,K(r)} (t)
=(1-t^3)(1-3t^3 + t^6)(1+t^3+ t^6)^2$.\\
(4) For  $r=29/75, 
\widetilde{\Delta}_{\rho_0,K(r)} (t)
=(1-t^3)(4-t^3)(1-4t^3)$.\\
(5) For $r=227/777$,
\begin{align*}
\widetilde{\Delta}_{\rho_0,K(r)} (t)
&=(1-t^3)(1-3t^3 +t^6)(1+t^3 +t^6)(2-3t^3 +2t^6)\\
&\ \ \times(4- 36t^3-35t^6
-71t^9 -35t^{12}-36t^{15}+4t^{18}).
\end{align*}

\end{ex}

%Sec 6 
\section{Metabelian representations}

In this section, 
we discuss general metabelian representations of knot groups.
First we define metabelian groups on which the knot groups map
\cite{Ha}.

Let $p$ be a prime and let $\Phi_n$ be the $n$-th cyclotomic polynomial. 
We assume that $(n,p)=1$ and $\Phi_n$ is irreducible over $\ZZ/p$. 
Let $k$ be the degree of $\Phi_n$.
Let $A(p,k)$ denote the elementary abelian $p$-group of order $p^k$, 
i.e. the direct product of $\ZZ/p$ $k$-times.
We define a semi-direct product 
$\ZZ/n \marubackslash A(p,k)$ on which the knot
groups map. For convenience, we denote 
$M(n|p,k) = \ZZ/n \marubackslash A(p,k)$.
Consider $A(p,k)$ as a $k$-dimensional vector space over $\ZZ/p$ and 
take a basis for $A(p,k)$, say $\{b_1, b_2, \cdots, b_k\}$. 
Let $T$ be the companion matrix of $\Phi_n$ over $\ZZ/p$.
Now, let $s$ be a generator of $\ZZ/n$ and fix it. 
An element $g$ of $M(n|p,k)$ is written as $g=s^{\ell} a_g$, 
where $a_g=b_1^{\lambda_1} b_2^{\lambda_2} \cdots b_k^{\lambda_k}$, 
$0 \leq {\ell}< n$ , and $0 \leq {\lambda_j} < p$, $1 \leq j \leq k$ .
We define the action of $\ZZ/n$ on $A(p,k)$ by conjugation:
\begin{equation} %6.1
s a_g s^{-1}=a_h,\ {\rm where\ }
a_h=b_1^{\mu_1} b_2^{\mu_2} \cdots b_k^{\mu_k}\ {\rm  is\
given\ by}
\end{equation}
%(6.1) $s a_g s^{-1}$ = $a_h$ ,
\begin{equation}%6.2
(\lambda_1, \lambda_2, \cdots, \lambda_k ) T
 =( \mu_1,\mu_2, \cdots, \mu_k).
 \end{equation}
To define the twisted Alexander polynomial 
$\widetilde{\Delta}_{\rho,K}(t)$ associated to 
a metabelian representation $\rho$, we need a
faithful representation of $M(n|p,k)$ in $GL(m,\ZZ)$ for some $m$.
Now let $N= \{1,s,s^2, \cdots, s^{n-1}\}$ be a subgroup of $M(n|p,k)$
generated by $s$ and let $\widehat{N} = \{N1=N, N2, \cdots, Np^k\}$ be the
set of all right cosets of $M(n|p,k)$ mod $N$. 
Using (6.1), we see that
the right multiplication of $g \in M(n|p,k)$ on $\widehat{N}$ induces
a permutation representation $\sigma$ of $M(n|p,k)$ in $S_{p^k}$ and hence
in $GL(p^k, \ZZ)$ via permutation matrices that is denoted by $\xi$ .

Suppose that there is a homomorphism $f$ from $G(K)$ onto $M(n|p,k)$ for
some $n,p$ and $k$. (It is known that if such a homomorphism exists, then
$p$ divides $\prod{\Delta_K(\omega_j)}$,
where the product runs over all primitive $n$-th roots $\omega_j$ of $1$.)
Then $f$ induces a representation $\rho =f \circ \sigma \circ \xi: 
G(K)\longrightarrow M(n|p,k) \longrightarrow S_{p^k} \longrightarrow 
GL(p^k ,\ZZ)$.

We may assume without loss of generality that 
for one meridian generator $x$,
\begin{equation}
f(x)=s.
\end{equation}
%(6.3) $f (x) = s$.

Let $\widetilde{\Delta}_{\rho, K}(t)$ be the twisted Alexander
polynomial of a knot $K$
associated to $\rho$. 
Then we propose the following conjecture:
%Conjecture 6.1. 
\begin{yosou}\label{conj:6.1}
$\widetilde{\Delta}_{\rho, K}(t)$ is of the form:
\begin{equation}%(6.4)
\widetilde{\Delta}_{\rho, K}(t) =
\left[\frac{\Delta_{K}(t)}{1-t}\right]
\varphi (t),\ 
{\rm where\ } \varphi (t)\ {\rm  is\ an\ integer\ polynomial\ in\ }
t^n.
\end{equation}
\end{yosou}

In the rest of this section, we discuss several examples that support this
conjecture.

%Example 6.2. 
\begin{ex}
Consider a torus knot $K(1/p)$, $p$ being an odd prime.
Let $G(K(1/p))=\langle x,y| Wx=yW\rangle$, $W=(xy)^{\frac{p-1}{2}}$, 
be a Wirtinger presentation.
Since $\Phi_p = 1+z+z^2+ \cdots +z^{p-1}$ is 
irreducible over $\ZZ/2$,
$G(K(1/p))$ maps on $M(p|2,p-1)$ by a mapping $f(x)=s$ and $f(y)=sb_1$,
and hence we obtain
a metabelian representation $\rho_p$ of $G(K(1/p))$,
%\begin{equation*}
$\rho_p : G(K(1/p)) \rightarrow M(p|2,p-1) \rightarrow
S_{2^{p-1}} \rightarrow GL(2^{p-1},\ZZ)$.
%\end{equation*}

For example, if $p=3$, then $M(3|2,2)$ is the alternating group $A_4$ and
$\rho_3$ coincides with $\rho$ defined in (1.2).
Let $\widetilde{\Delta}_{\rho_3, K(1/p)}(t)$ be the twisted Alexander
polynomial of $K(1/p)$ associated to $\rho_3$. Then computations show that

%(6.5)
\begin{align} 
&(1) \widetilde{\Delta}_{\rho_3, K(1/3)}(t)=
\left[\frac{\Delta_{K(1/3)}(t)}{1-t}\right] (1-t^3).\nonumber\\
&(2) \widetilde{\Delta}_{\rho_5, K(1/5)}(t)=
\left[\frac{\Delta_{K(1/5)}(t)}{1-t}\right] (1-t^5)^5 (1+t^5)^4.
\end{align}

It is quite likely that for any odd prime $p$,
%(6.6) 
\begin{equation}
\widetilde{\Delta}_{\rho_p, K(1/p)}(t) =
\left[\frac{\Delta_{K(1/p)}(t)}{1-t}\right] (1-t^p)^m (1+t^p)^{m-1},\
{\rm where\ }m=2^{p-2}-\left[\frac{2^{p-1}-1}{p}\right].
\end{equation}
\end{ex}

%Example 6.3. 
\begin{ex}
Consider $M(4|3,2)=Z/4 \marubackslash (\ZZ/3 \oplus \ZZ/3)$.
Since $\Phi_4 = 1+z^2$ is irreducible over $\ZZ/3$, the group has the
following presentation:
%(6.7) 
\begin{equation}
M(4|3,2)= \langle s, a, b| s^4 = a^3 = b^3 =1, ab=ba, sa s^{-1}
= b^{-1}, sb s^{-1} = a\rangle.
\end{equation}

Let $N = \{1,s , s^2, s^3\}$ 
be a subgroup of $M(4|3,2)$ generated by $s$.
Let $\widehat N = \{N=N1, N2,
\cdots, N9\}$ be the set of all right cosets of $M(4|3,2)$ mod $N$. Let
$\sigma$ be the
permutation representation of $M(4|3,2)$ in $S_9$. 
For example,
$\sigma(s)= (1)(2435)(6897)$ and $\sigma(sa)=(1264)(5378)(9)$.

Now let $G(K(r)) =\langle x,y| R_r = 1\rangle$ be a Wirtinger presentation of
$G(K(r))$. Then the groups of $2$-bridge knots $K(3/5), K(3/7), K(5/13),
K(11/17)$ and $K(13/23)$ map onto $M(4|3,2)$ by $f(x) =s$ and $f(y)=sa$
and
a homomorphism : $\rho(x) = \xi(\sigma(s))$ and $\rho(y)=
\xi(\sigma(sa))$ defines a metabelian representation of $G(K(r))$. The
twisted Alexander
polynomials $\widetilde{\Delta}_{\rho, K(r)}(t)$ associated to $\rho$ are
of the form : $\left[\frac{\Delta_{K(r)}(t)}{1-t}\right] \varphi_{K(r)} (t)$, and
$\varphi_{K(r)} (t)$ is given as follows.
%(6.8)
\begin{align}
&(1)\ \varphi_{K(3/5)} (t) = (1-t^4)^2,\nonumber\\
&(2)\ \varphi_{K(3/7)} (t)=4(1-t^4)^2, \nonumber\\
&(3)\ \varphi_{K(5/13)} (t)= (1-t^{12})^2,\nonumber\\
&(4)\ \varphi_{K(11/17)} (t) = (1-t^4)^4 (1+t^4+t^8)^3\ {\rm and}\nonumber\\
&(5)\ \varphi_{K(13/23)} (t) = (1-t^4)^2
(4-13t^4-9t^8-13t^{12}+4t^{16}).
\end{align}
\end{ex}

%Example 6.4. 
\begin{ex}
Consider $M(3|5,2) = \ZZ/3 \marubackslash (\ZZ/5 \oplus \ZZ/5)$.
Since $\Phi_3=1+z+z^2$ is irreducible over $\ZZ/5$, 
the group has a presentation:
$\langle s, a, b| s^3 = a^5 = b^5 =1, ab=ba, sa s^{-1} = b^{-1}, sb
s^{-1}=a b^{-1}\rangle$.
\end{ex}

As is shown in Example 6.3, $M(3|5,2)$ can be represented in $S_{25}
\subset GL(25,\ZZ)$.
The groups of $2$-bridge knots $K(3/7), K(7/11), K(9/23)$ and $K(9/31)$
map onto $M(3|5,2)$ by $f(x) = s$ and $f(y) = sa$. The twisted
Alexander polynomials associated to $\rho=f \circ \sigma \circ \xi$ :
$G(K(r)) \longrightarrow GL(25,\ZZ)$ are of the form:
$\left[\frac{\Delta_{K(r)}(t)}{1-t}\right] \varphi_{K(r)} (t)$ and we have
%(6.9)
\begin{align}
&(1)\ \varphi_{K(3/7)} (t) = 16(1-t^3)^8,\nonumber\\
&(2)\ \varphi_{K(7/11} (t) = (1-t^3)^8
(1-3t^3-2t^6-6t^9-5t^{12}-6t^{15}-2t^{18}-3t^{21}+t^{24})^2,\nonumber\\
&(3)\ \varphi_{K(9/23)} (t) = \nonumber\\
&\ \ (1-t^3)^8 (1-5t^6-20t^9
-28t^{12}-20t^{15}-5t^{18}+t^{24})^2\nonumber\\
&\hspace*{1.5cm}\times(1-5t^3+10t^6-10t^9+7t^{12}-10t^{15}+10t^{18}-5t^{21}+t^{24})^2\ {\rm and}\nonumber\\
&(4)\ \varphi_{K(9/31)} (t) = \nonumber\\
&\ \ (1-t^3)^8
(1+3t^3-6t^6+15t^9-15t^{12}+15t^{15}-6t^{18}+3t^{21}+t^{24})^2\nonumber\\
&\ \ \times
(4+12t^3+36t^6+30t^9+35t^{12}+30t^{15}+36t^{18}+12t^{21}+4t^{24})^2.
\end{align}

%Example 6.5. 
\begin{ex}
Consider $M(4|5,2) = \langle s, a, b| s^4 = a^5 = b^5 =1,
ab=ba, sa s^{-1} = b^{-1}, sb s^{-1} =a\rangle$. 
Since $\Phi_4$ is reducible over $\ZZ/5$, 
we cannot apply the previous argument,
but we see that $M(4|5,2)$ can be represented in $S_{25} \subset
GL(25,\ZZ)$ and the group of a $2$-bridge knot $K(5/9)$ maps onto 
$M(4|5,2)$ by $f(x)=s$ and $f(y)=sa$ and hence $G(K(5/9))$ 
has a metabelian representation $\rho$ in $GL(25,\ZZ)$.
The twisted Alexander polynomial of $K(5/9)$ associated to $\rho$ is
given by
%(6.10) 
\begin{equation}
\widetilde{\Delta}_{\rho, K(5/9)}(t)=
\left[\frac{\Delta_{K(5/9)}(t)}{1-t}\right] 16(1-t^4)^6.
\end{equation}
\end{ex}

In the last three examples below, we consider a non-rational knot $K =8_5$
and non-alternating knots $10_{145}$ and $10_{159}$ in Rolfsen table,
and show that Conjecture \ref{conj:6.1} holds for $K$.

%%%%\begin{ex}
%%%%Let $K$ be the knot $8_5$. The group $G(K)$ has a Wirtinger presentation
%%%%\cite[Example 10.5]{HM2}:
%%%%$G(K)=\langle x,y,z|R_1,R_2\rangle$, where
%%%%$R_1 =(x^{-1}y^{-1}zyxy^{-1}x^{-1}y^{-1})x(yxyx^{-1}y^{-1}z^{-1}yx)y^{-1}$, 
%%%%and
%%%%$R_2 = (xzyxy^{-1})z(yx^{-1}y^{-1}z^{-1}x^{-1})y^{-1}$.
%%%%\end{ex}

%Example 6.6. 

\begin{ex}
Let $K$ be the knot $8_5$. 
The group $G(K)$ has a Wirtinger presentation 
\cite[Example 10.5]{HM2}: 
\begin{align*}
&G(K)=\langle x,y,z|R_1,R_2\rangle,\ {\rm  where}\\
&R_1 =(x^{-1}y^{-1}zyxy^{-1}x^{-1}y^{-1})x(yxyx^{-1}y^{-1}z^{-1}yx)y^{-1},\ {\rm  and}\\
&R_2 = (xzyxy^{-1})z(yx^{-1}y^{-1}z^{-1}x^{-1})y^{-1}.
\end{align*}
We see that $G(K)$ maps on a metabelian group $M(3|2,2)$ that is an
alternating group $A_4$. (We note that $\Delta_K(t)
=(1-t+t^2)(1-2t+t^2-2t^3+t^4))$.
A mapping $f(x)=f(z)=(123)$ and $f(y)=(142)$ induces an
$A_4$-representation $\rho$ of $G(K)$ into $GL(4,\ZZ)$. 
The twisted Alexander polynomial is given by
%(6.11)
\begin{equation} 
\widetilde{\Delta}_{\rho, K}(t)
=\left[\frac{\Delta_K(t)}{1-t}\right](1-t^3)(1-8t^3-6t^6-8t^9+t^{12}).
\end{equation}
\end{ex}

%Example 6.7. 
\begin{ex}
Let $K$ be the knot $10_{145}$. The Alexander polynomial
$\Delta_K(t)$ is $1~+~t-3t^2+t^3+t^4$.
The group $G(K)$ has a Wirtinger presentation:
\begin{align*}
&G(K)=\langle x,y,z|R_1,R_2\rangle,\ {\rm where}\\
&R_1 = (y^{-1}z x^{-1}z^{-1}y z x
y^{-1}x^{-1}z^{-1})y(zxyx^{-1}z^{-1}y^{-1}zxz^{-1}y)z^{-1},\ {\rm  and}\\
&R_2 =(z^{-1}y^{-1}zx^{-1}z^{-1}yzxy^{-1})z(yx^{-1}z^{-1}y^{-1}zxz^{-1}yz)x^{-1}.
\end{align*}
Now we see that $G(K)$ maps on the group $M(5|2,4) = \ZZ/5 
\marubackslash
(\ZZ/2)^4$. This group has a presentation:
%(6.12) 
\begin{align}
\langle s,b_1,b_2,b_3,b_4|
&s^5=b_i^2 = 1, b_i b_j=b_j b_i,1 \leq i, j \leq 4, \nonumber\\
&s b_1 s^{-1} = b_4^{-1}, s b_2 s^{-1} = b_1 b_4^{-1},\nonumber\\
&s b_3 s^{-1} = b_2 b_4^{-1},s b_4 s^{-1} = b_3b_4^{-1}\rangle.
\end{align}
We see easily that a mapping $f(x)=sb_1b_2 b_3 b_4$, $f(y)=s b_1$ and
$f(z) = s$ is, in fact, a homomorphism and we can represent $M(5|2,4)$ in
$GL(16,\ZZ)$. Then the twisted Alexander polynomial of $K$
associated to a representation $\rho : G(K) \longrightarrow GL(16,\ZZ)$ is
given as follows.
%(6.13) 
\begin{equation}
\widetilde{\Delta}_{\rho, K}(t) = \left[\frac{\Delta_K(t)}{1-t}\right]
(1-t^5) (1+14t^5 + t^{10}) (1+3t^5+t^{10}).
\end{equation}
\end{ex}

%Example 6.8. 
\begin{ex}
Let $K$ be the knot $10_{159}$. The Alexander polynomial
$\Delta_K(t)$ is $(1-t+t^2)(1-3t+5t^2-3t^3+t^4)$.
The group $G(K)$ has a Wirtinger presentation:
$G(K)=\langle x,y,z|R_1,R_2\rangle$, where
$R_1 = (xzx^{-1}z^{-1}y^{-1}zy)x(y^{-1}z^{-1}yzxz^{-1}x^{-1})y^{-1}$,
and
$R_2 =
(x^{-1}yxzx^{-1}z^{-1}x^{-1}y^{-1}z^{-1})y(zyxzxz^{-1}x^{-1}y^{-1}x)z^{-1}$.
We see that $G(K)$ maps on two metabelian groups $M(3|2,2)$ and
$M(5|2,4)$. 
The first group is an alternating group $A_4$ and a
mapping $f(x)=f(y)=(123)$ and $f(z)=(142)$ induces an
$A_4$-representation $\rho$ of $G(K)$ into $GL(4,\ZZ)$. 
The twisted Alexander polynomial is given by
%(6.14) 
\begin{equation}
\widetilde{\Delta}_{\rho, K}(t) = 
\left[\frac{\Delta_K(t)}{1-t}\right]
(1-t^3)(1-3t^3-3t^6-3t^9+t^{12}).
\end{equation}
Now consider the second group $M(5|2,4) =\ZZ/5 
\marubackslash (\ZZ/2)^4$. 
This group has a presentation (6.12).
We see easily that a mapping $f(x)=s, f(y)=s b_1b_4$ and $f(z)=s b_1$ is a
homomorphism. 
As before, we represent $M(5|2,4)$ in $GL(16,\ZZ)$. 
Then the twisted Alexander polynomial of $K$
associated to a representation $\rho : G(K) \longrightarrow GL(16,\ZZ)$ is
given as follows.
%(6.15) 
\begin{equation}
\widetilde{\Delta}_{\rho, K}(t) = \left[\frac{\Delta_K(t)}{1-t}\right]
\varphi(t),\ {\rm where\ } 
\end{equation}
$\varphi(t) =
(1-t^5)(1+3t^5+t^{10})(1-31t^5+12t^{10}-31t^{15}+t^{20})(1+5t^5+52t^{10}+5t^{15}+t^{20}).
$
\end{ex}

\medskip
\noindent
{\bf Acknowledgements. } 
The first author is
partially supported by MEXT, Grant-in-Aid for
Young Scientists (B) 18740035,
and the second author is
partially supported by NSERC Grant~A~4034.

\end{document}